\documentclass[11pt,a4paper]{article}
\usepackage{dsfont}
\usepackage{parskip}
\usepackage{amssymb}
\begin{document}
\author{I.R. Kayumov, K.-J. Wirths}
\title{On the sum of squares of the coefficients of  Bloch functions}
\makeatletter
\renewcommand{\@makefnmark}{}
\makeatother
\newtheorem{lemma}{Lemma}
\newtheorem{remark}{Remark}
\newtheorem{corollary}{Corollary}
\newtheorem{theorem}{Theorem}
\newtheorem{proposition}{Proposition}
\newtheorem{df}{Definition}

\maketitle
\begin{footnote}{Keywords: Bloch functions, Subordination, Area functional.\\ AMS classification numbers: 30H30, 30B10.}
\end{footnote}
\makeatletter
\renewcommand{\@makefnmark}{}

\makeatother

\maketitle
\begin{abstract}
In this article several types of inequalities for weighted sums of the moduli of Taylor coefficients for Bloch functions are proved.
\end{abstract}

\section{Introduction.}
It was one of the most famous items in function theory, when A.Bloch proved in \cite{B} the following theorem\\
{\bf Theorem A.} {\it Let} $F$ {\it be analytic in the unit disc} $\mathds{D}$ {\it and} $F'(0)\,=\,1.$ {\it Let further} $B_F$ {\it be the supremum of all numbers} $r>0$ {\it with the following property. There exist a complex number} $a$ {\it and a domain} $\Omega \subset \mathds{D}$ {\it such that} $F$ {\it is injective on} $\Omega$ {\it and} $F(\Omega)\,=\,\{z\,\mid\, |z-a|<r\}$. {\it The infimum of these numbers} $B_F$ {\it is positive.}

\bigskip
\noindent
Since then, this infimum has been called Bloch's constant $B$. In other words, the range of any function of the above type covers a schlicht disc with radius $B$.  The exact value of Bloch's constant is not known.\\
In \cite{L}, E. Landau proved that it suffices for estimations of $B$ to  consider only those functions that satisfy in addition to the above conditions the inequality
\begin{equation} \label{Bloch}
 |F'(z)| \leq \frac{1}{1-|z|^2}, \quad z\in \mathds{D}.
 \end{equation}
This was the beginning of research on the class of Bloch functions that are analytic in $\mathds{D}$ and satisfy the condition
\[ \sup_{z\in \mathds{D}}  (1-|z|^2)|F'(z)|\,<\,\infty. \]
This class becomes a normed vector space, the so called Bloch space, if it is endowed with the Bloch norm
\[
\|F\|\,:=\,|F(0)|\,+\, \sup_{z\in \mathds{D}}  (1-|z|^2)|F'(z)|.
\]
Ch. Pommerenke and his co-authors proved  a number of important theorems on the Bloch space (compare for example \cite{AnClPo} and \cite{Po}), among them some theorems on the behaviour of the coefficients $b_k, k\to \infty,$ in the expansions
\[
F(z)\,=\, \sum_{k=0}^{\infty}b_kz^k,\]
where $F$ is in the Bloch space.\\
Another area of research was related more closely to the condition (\ref{Bloch}). It concerned the class
\[
\mathcal{B}\,:=\,\left\{F\,\mid \, |F'(z)| \leq \frac{1}{1-|z|^2},\, z\in \mathds{D}\right\}\]
and the behaviour of its elements. In \cite{R} the following problem was posed\\
{\bf Problem N:} (see \cite{R}) {\it If} $f$ {\it is analytic in the unit disc} $\mathds{D}$ {\it and}
\[
\left|f\left(re^{i\theta}\right)\right|\,<\,\frac{1}{1-r^2},\]
{\it for what positive real numbers} $\alpha,$ {\it is it true that}
\[\left|f'\left(re^{i\theta}\right)\right|\,<\,\frac{\alpha}{(1-r^2)^2}\, ?\]
{\it It is known that} $\alpha = 4 $ {\it will suffice.}\\
{\it What is the coefficient region for this class ?}\\
The first part of this problem was solved in \cite{Wir} and \cite{AvkKay} with different methods. Partial answers to the second question can be found in \cite{Wir}, \cite{bonk}, \cite{CG},  \cite{Ro}, and \cite {ST}. In the following we shall use \\
{\bf Theorem B.} (see \cite{Wir} and \cite{bonk}) {\it Let} $F\in \mathcal{B}$ {\it and} $x\in [0,1/\sqrt{3}]$ {\it such that}
\[
|b_1|\,= \,\frac{3\sqrt{3}}{2}x(1\,-\,x^2).\]
{\it Then this equation and the inequality}
\[
|b_2|\,\leq \,\frac{3\sqrt{3}}{4}(1\,-\,3x^2)(1\,-\,x^2) \]
{\it describe the coefficient region} $\{(b_1,b_2)\,\mid\,F\in \mathcal{B}\}.$

\bigskip
\noindent
In \cite{bonk} and \cite{CG} an intimate relation between knowledges on this coefficient region of $\mathcal{B}$ and calculations of estimates for Bloch's constant is revealed.

\bigskip
\noindent
Another way to get information on the coefficient region of $\mathcal{B}$ consists in the consideration of weighted sums of moduli of Taylor coefficients. One example for an estimate of this type is the use of Parseval's formula (see f. i. \cite{Po})
\begin{equation}\label{Basic}
\sum_{k=1}^{\infty}k^2|b_k|^2r^{2k-2}\,=\,\frac{1}{2\pi}\int_0^{2\pi}\left|F'\left(re^{i\theta}\right)\right|^2\,d\theta\,\leq \,\frac{1}{(1-r^2)^2},\,\,r\in (0,1).
\end{equation}
Another inequality of this type can easily proved using the maximum principle (compare \cite{RW}).\\
{\bf Proposition 1.} {\it Let} $F\in \mathcal{B}$. {\it Then}
\[
\sum_{k=n+1}^{\infty}k^2|b_k|^2r^{2k-2}\,\leq \,\frac{(n+2)^{n+2}}{4n^n}r^{2n},\quad 0\leq r \leq \sqrt{\frac{n}{n+2}}\,=\,r_n.\]
{\it In this estimate, equality is attained for}
\[
F'_n(z)\,=\,\frac{n+2}{2}\left(\frac{n+2}{n}\right)^{\frac{n}{2}}z^n.\]
{\it Proof.} Let us set
$$
\Phi_n(z)=\sum_{k=n+1}^{\infty}k^2|b_k|^2z^{2k-2}.
$$ From (\ref{Basic}) it follows that
$$
|\Phi_n(z)| \leq \frac{1}{(1-r^2)^2} = \frac{(n+2)^{n+2}}{4n^n}r^{2n}, \quad |z|=r=r_n.
$$ Consequently, by the maximum principle
$$
\left|\frac{\Phi_n(z)}{z^{2n}}\right| \leq  \frac{(n+2)^{n+2}}{4n^n}, \quad |z| \leq r_n.
$$ Setting $z=r \leq r_n$ we prove Proposition 1.

\bigskip
\noindent
{\bf Remark 1.} If Proposition 1 is compared with (\ref{Basic}), it is obvious that in (\ref{Basic}) equality is attained for $r\,=\,r_n$ and $F\,=\,F_n$. It is an open question whether there exist other values of $r$ such that in (\ref{Basic}) occurs equality.

\bigskip
\noindent
Inequalities of different type for weighted sums can be found in  \cite{K}.\\
In the paper \cite{IK}, another application of such sums was revealed. Namely, it was shown that the function
$$
\Phi(r):=\sup_{F \in \cal{B}}(1-r^2)\sum_{k=1}^\infty k|b_k|²r^{2(k-1)}
$$ is a decreasing function in $r$ and $\Phi(0.4)<0.9$. Moreover, in the cited paper it was demonstrated that
$$
\sigma^2_F \leq \Phi(r), \quad r \in [0,1),
$$ where
\[\sigma_F^2\,=\,\limsup_{r \to 1}\frac{1}{2\pi|\log(1-r)|)}\int_{|z|=r}|F(z)|^2|dz| \]
is the asymptotic variance of the Bloch function $F$. The asymptotic variance plays an important role in the probabilistic behaviour
of the Bloch function $F$. Namely, Makarov's law of the iterated logarithm (see \cite{makarov85}, and also \cite{IK}, \cite{kayumov}) asserts that
$$
\limsup_{r \to 1} \frac{ |F'(re^{i\theta})|}{\sqrt{\log \frac{1}{1-r} \log \log\log  \frac{1}{1-r}}}  \leq \sup_{F \in \cal{B}}\sqrt{\sigma^2_F} \quad \mbox{for almost all } \theta \in [0,2\pi).
$$
The present paper is dedicated to further estimates that generalize and sharpen some of the above ones.

\bigskip
\noindent

\section{Statement and proofs of the results.}

 \begin{theorem} \label{Th1} \,\,{\it Let} $F(z)\,=\,\sum_{k=0}^{\infty} b_kz^k\in \mathcal{B}, b_1=a\in (0,1),$
\[
a\,=\,\frac{3\sqrt{3}}{2}x(1-x^2),\quad x\in (0,1/\sqrt{3}),\]
\begin{equation} \label{Additional1}
0\,<\,r\,\leq \frac{\sqrt{1/3}-x}{1-x\sqrt{1/3}}.
\end{equation}

{\it Then the  sharp inequalities}
\[\sum_{k=1}^\infty k^2|b_k|^2 r^{2k} \leq \]

\begin{equation}\label{B}
\frac{27 r^2  (1 - x^2)^2 ((r^2+x^2)(1-r^2x^2)^2-6r^2x^2(1-x^2)(1-r^2))}{
 4 (1 - r^2 x^2)^5}, \end{equation}
{\it and}
\begin{equation}\label{B2}
\sum_{k=1}^\infty k|b_k|^2 r^{2k}\,\leq \frac{27 r^2 (1 - x^2)^2 (3x^2(1-r^2)^2+(1-x^2r^2)(r^2-x^2))}{8 (1- r^2 x^2)^4}
\end{equation}
{\it are valid.}
\end{theorem}

\bigskip
\noindent
{\it Proof.}\,\, Using the rotation $e^{-i\theta}F(e^{i\theta}z)$ it is easy to check that
\[
\{F'(z)\mid |z|=r, F\,\,{\rm{ as\,\, above}}\}\,=\,\{F'(r)\mid F\,\,{\rm{as\,\, above}}\}.\]
Following \cite{bonk}, Satz 2.2.1, this set is given by the bounded region whose boundary is the Jordan curve
\[
\{G'(re^{i\phi})\mid \phi \in [0,2\pi]\}\]
where
\[
G'(z)\,=\,-\frac{a}{x}\frac{z-x}{(1-zx)^3}\, = \sum_{k=1}^{\infty} kA_k z^{k-1}\]\[=\, \frac{a}{2}\sum_{k=1}^{\infty} k x^{k-3} (2 x^2 + (k-1) (x^2-1))z^{k-1}.\]

From   \cite{bonk}, Lemma 2.3.3 and Lemma 2.3.5 it follows that $G'(rz)$ maps the closed unit disc ${\overline{\mathds{D}}}$ univalently. Since $F'(0)=G'(0)$ and $F'(r{\overline{\mathds{D}}})\subset G'(r{\overline{\mathds{D}}}),$ this means that $F'(rz)$ is subordinate to $G'(rz)$. According to a famous Lemma of Rogosinski (see \cite{Rog}), this implies
\[\sum_{k=1}^nk^2|b_k|^2²\rho^{2(k-1)}\,\leq \sum_{k=1}^nk^2|A_k|^2²\rho^{2(k-1)},\,\,\rho\in [0,r].\]
The limiting process $n\to \infty$ and a straightforward calculation delivers (\ref{B}).\\
If we let $u=\rho^2$ and integrate this inequality with respect to $u$ from $0$ to $r^2$, we get the assertion (\ref{B2}).

\bigskip
\noindent
{\bf Remark 2.} If we assume $F'(0)=\epsilon a,\,\,|\epsilon|=1$, we may apply the above reasoning to $\overline{\epsilon}F(z)$ and  $\overline{\epsilon}F_x(z)$ and we see that (\ref{B}) and (\ref{B2}) are valid likewise.

\bigskip
\noindent
{\bf Remark 3.} Another subordination theorem for Bloch functions may be found in \cite{BMY}.

\begin{theorem} \label{Th2} Let $F\in \mathcal{B}$. Then the  sharp inequality
$$
\sum_{k=1}^\infty k^2|b_k|^2 r^{2k} \leq \frac{27}{4} r^4, \quad \sqrt{4/15} \leq r \leq \sqrt{1/3},
$$
is valid.
\end{theorem}
Proof. If we let $r=\frac{1}{\sqrt{3}}$ in (\ref{Basic}), we get
$$
\sum_{k=3}^\infty k^2|b_k|^2 \frac{1}{3^{k-1}} \leq \frac{9}{4} - |b_1|^2-(4/3)|b_2|^2.
$$
According to Remark 1, this inequality is sharp. From here we see that
$$
\sum_{k=3}^\infty k^2|b_k|^2r^{2(k-1)} \leq 9r^4 \left(\frac{9}{4} - |b_1|^2-(4/3)|b_2|^2 \right).
$$ and consequently
$$
\sum_{k=1}^\infty k^2|b_k|^2r^{2(k-1)} \leq |b_1|^2+4|b_2|^2r^2+ 9r^4 \left(\frac{9}{4} - |b_1|^2-(4/3)|b_2|^2 \right)
$$
\[
=(1-9r^4)|b_1|^2+(4r^2-12r^4)|b_2|^2 + \frac{81r^4}{4} \leq \frac{27r^2}{4}, \quad r \ge \sqrt{4/15}.
\]
From Theorem B we know that it is sufficient to verify this inequality in the case when
$$
b_1\,=\,\frac{3\sqrt{3}}{2}x(1-x^2),\quad x\in (0,1/\sqrt{3}),
$$ and
$$
b_2=\frac{3\sqrt{3}}{4} (1 - x^2) (1 - 3 x^2).
$$

We have
$$
(1-9r^4)|b_1|^2+(4r^2-12r^4)|b_2|^2 + \frac{81r^4}{4} - \frac{27r^2}{4}=
$$
$$
(27/4) (1 - 3 r^2) x^2 (1 - 2 x^2 + x^4 +
   r^2 (-5 + 16 x^2 - 21 x^4 + 9 x^6)).
$$ It remains to show that
$$
1 - 2 x^2 + x^4 +
   r^2 (-5 + 16 x^2 - 21 x^4 + 9 x^6) \leq 0, \quad  1/\sqrt{3} \ge r \ge \sqrt{4/15}.
$$ Obviously,  the first term $1-2x^2+x^4$ is positive, whereas the second one is negative for $x\in (0,1/\sqrt{3})$. Hence, it is evidently enough to verify this inequality at $ r= \sqrt{4/15}$. In this case   we have
$$
1 - 2 x^2 + x^4 +
   r^2 (-5 + 16 x^2 - 21 x^4 + 9 x^6) = (1/15) (-1 + 3 x^2)^2 (-5 + 4 x^2) \leq 0.
$$ Theorem \ref{Th2} is proved.

\bigskip
\noindent
{\bf Problem 1.} Which is the biggest interval $\left[c_1, \frac{1}{\sqrt{3}}\right], c_1>0,$ such that the inequality of Theorem 2 remains valid in this interval ?

\bigskip
\noindent
If one adds to  the inequality (\ref{Bloch}) the slightly stronger inequality
\begin{equation}\label{B4}
|b_1|+\left|\sum_{k=2}^{\infty} kb_k  z^{k-1}\right| \leq \frac{1}{1-|z|^2}, \quad |z|=1/\sqrt{3},
\end{equation}
it is possible to improve the length of the interval from Theorem 2.

\begin{theorem} \label{Th3} Suppose that a function $F\in\mathcal{B}$ satisfies (\ref{B4}). Then the  sharp inequality
$$
\sum_{k=1}^\infty k^2|b_k|^2 r^{2k} \leq \frac{27}{4} r^4, \quad \sqrt{(9-\sqrt{65})/6} \leq r \leq \sqrt{1/3},
$$
is valid.
\end{theorem}
Proof. We have
$$
\sum_{k=2}^\infty k^2|b_k|^2 \frac{1}{3^{k-1}} \leq \left(\frac{3}{2} - |b_1|\right)^2
$$ and hence
$$
\sum_{k=3}^\infty k^2|b_k|^2 \frac{1}{3^{k-1}} \leq \left(\frac{3}{2} - |b_1|\right)^2-(4/3)|b_2|^2.
$$
From here we see that
$$
\sum_{k=3}^\infty k^2|b_k|^2r^{2(k-1)} \leq 9r^4 \left[ \left(\frac{3}{2} - |b_1|\right)^2-(4/3)|b_2|^2\right].
$$ and consequently
$$
\sum_{k=1}^\infty k^2|b_k|^2r^{2(k-1)} \leq |b_1|^2+4|b_2|^2r^2+ 9r^4 \left[ \left(\frac{3}{2} - |b_1|\right)^2-(4/3)|b_2|^2\right] \leq \frac{27r^2}{4}.
$$
In view of Theorem B it is enough to verify this inequality in the case when
$$
b_1\,=\,\frac{3\sqrt{3}}{2}x(1-x^2),\quad x\in (0,1/\sqrt{3}),
$$ and
$$
b_2=\frac{3\sqrt{3}}{4} (1 - x^2) (1 - 3 x^2).
$$ Also we may suppose that
$r=\sqrt{(9-\sqrt{65})/6}$. We have
$$ |b_1|^2+4|b_2|^2r^2+ 9r^4 \left[ \left(\frac{3}{2} - |b_1|\right)^2-(4/3)|b_2|^2\right] - \frac{27r^2}{4}=$$
$$
-\frac{9x}{8}(2 \sqrt{3} (73 - 9 \sqrt{65}) + (-737 + 91 \sqrt{65}) x +
   2 \sqrt{3} (-73 + 9 \sqrt{65}) x^2 +
   $$
   $$
   + (1858 - 230 \sqrt{65}) x^3 +
   3 (-587 + 73 \sqrt{65}) x^5 - 72 (-8 + \sqrt{65}) x^7)=
$$
$$
(1-\sqrt{3}x)^2\left(\frac{9}{4} \sqrt{3} (-73 + 9 \sqrt{65} ) x + \frac{9}{8} (-139 + 17 \sqrt{65} ) x^2 -
 \frac{9}{8} \sqrt{3}  (-153 + 19 \sqrt{65} ) x^3 - \right.
 $$
 $$ \left.
- \frac{9}{8} (-395 + 49 \sqrt{65} ) x^4 + 18 \sqrt{3}  (-8 + \sqrt{65} ) x^5 +
 27 (-8 + \sqrt{65} ) x^6 \right)=
$$
$$
(1-\sqrt{3}x)^2(-1.71348 \ldots x - 2.18432\ldots x^2 - 0.712771\ldots x^3 -
$$
$$
 -0.0569584 \ldots x^4 + 1.941\ldots x^5 + 1.68096 \ldots x^6).
$$
From here we see that negative coefficients are dominating. Consequently, this polynomial is negative in the interval $(0,1/\sqrt{3})$.  Theorem \ref{Th3} is proved.

\bigskip
\noindent
{\bf Problem 2.}  Which is the biggest interval $\left[c_2, \frac{1}{\sqrt{3}}\right], c_2>0,$ such that the inequality of Theorem 3 remains valid in this interval ?

Now we are going to study the behavior of the area functional
$$
\sum_{k=1}^\infty k|b_k|^2 ²r^{2k}.
$$ A simple integration of the inequality (\ref{Basic}) gives us the inequality
$$
\sum_{k=1}^\infty k|b_k|^2 ²r^{2k} \leq \frac{r^2}{1-r^2}
$$ which is, unfortunately, not sharp for all $r \in (0,1)$.

\begin{theorem} \label{Th4} \,\,Let $F\in\mathcal{B}$ be as in Theorem 1 and let
$r\leq \frac{1}{\sqrt{3}}.$\\

Then for any $n\in \mathds{N}\setminus\{1\}$ the inequality
\begin{equation}\label{Bl}
\sum_{k=2}^nk|b_k|^2 ²r^{2k}\,\leq \frac{3(9-4a^2)^2}{64a^4}\sum_{k=2}^n\frac{1}{k}\left(\frac{4}{3}a^2r^2\right)^k\end{equation}
is valid.
\end{theorem}

\noindent
{\it Proof.}\,\,Since $|F'(z)|\leq \frac{1}{1-|z|^2}\leq \frac{3}{2}$ for $r\leq\frac{1}{\sqrt{3}},$ the image of the unit disc under $F'(z/\sqrt{3})$ lies in the disc with radius $3/2$ around the origin. Hence, this function is subordinate to the function
\[
H(z)\,=\,\frac{3}{2}\frac{z\,+\,\frac{2a}{3}}{1\,+\,\frac{2a}{3}z}.\]
As
\[
F'\left(\frac{z}{\sqrt{3}}\right)\,=\,a\,+\,\sum_{k=2}^{\infty}k\,b_k\left(\frac{z}{\sqrt{3}}\right)^{k-1}\]
and
\[
H(z)\,=\,a\,+\,\sum_{k=2}^{\infty}\frac{3}{2}\frac{9-4a^2}{9}\left(-\frac{2a}{3}\right)^{k-2}z^{k-1},\]
we get from Rogosinski's theorem (see again \cite{Rog})
\[
\sum_{k=2}^{n}k^2|b_k|^2\frac{1}{3^{k-1}}\,\leq \frac{9(9-4a^2)^2}{64a^4}\sum_{k=2}^{n}\left(\frac{4a^2}{9}\right)^k.\]
Now we   summands of both sums by $\lambda_k\,=\,\frac{(3r^2)^k}{k}$. Since this is a  decreasing sequence of nonegative numbers, using properties of the Abel transformation (see also \cite{DurenUniv-83-8}, Theorem 6.3), we get
\[3\sum_{k=2}^{n}k|b_k|^2r^{2k}\leq  \frac{9(9-4a^2)^2 }{64a^4}\sum_{k=2}^{n}\left(\frac{4a^2}{3}\right)^k\frac{r^{2k}}{k}.\]
This results in the inequality (\ref{Bl}).

\bigskip
\noindent
{\bf Corollary 1.}\,\, {\it The limiting process} $n\to \infty$ {\it in this theorem results in the inequality}
\[\sum_{k=2}^{\infty}k|b_k|^2r^{2k}\leq \frac{3(9-4a^2)^2 }{64a^4}\left(-\log\left(1-\frac{4a^2r^2}{3}\right)\,-\,\frac{4a^2r^2}{3}\right)\,=\,B_a(r),\]
{\it where} $ r\leq\frac{1}{\sqrt{3}}.$

\bigskip
\noindent
An immediate consequence of Corollary 1 is \\
{\bf Corollary 2.}\,\, {\it For} $F\in \mathcal{B},\,\, r\leq \frac{1}{\sqrt{3}},$ {\it the inequality}
\[\sum_{k=2}^{\infty}k|b_k|^2r^{2k}\leq \frac{27 r^4}{8}\]
{\it is valid.}

\bigskip
\noindent
{\it Proof.}\,\,We have to prove that
\[
B_a(r)\,-\,\frac{27 r^4}{8}\leq 0,\quad r\leq \frac{1}{\sqrt{3}}.\]
To abbreviate the calculations, we let $w=\frac{4a^2r^2}{3}$. With this abbreviation the inequality to prove is the following
\[
H_a(w)\,=\,\left(1\,-\,\frac{4a^2}{9}\right)^2(-\log(1-w)\,-\,w)\,-\,\frac{w^2}{2}\,\leq \,0,\,\,w\in \left[0,\frac{4a^2}{9}\right].\]
Since $H_a(0)=0, H'_a(0)=0, H'_a(w)<0$ for sufficiently small positive values of $w$, and $H'$ has only one positive zero, it is sufficient to prove the above inequality for $w\,=\,\frac{4a^2}{9}$. If we let $v\,=\,\frac{4a^2}{9}$, this task reduces to the proof of the inequality
\[
-\log(1-v)\,-\,v\,-\frac{v^2}{2(1-v)^2}\,\leq \,0,\,\,v\in\left[0,\frac{4}{9}\right].\]
This proof can be done by elementary calculations.

{\bf Remark 4:} Corollary 2 follows immediately from the case $n=1$ of Proposition 1 using $k \leq k^2/2$ for $k\geq 2$.

\begin{theorem} \label{Th5} Let $F\in \mathcal{B}$. For
$$
R=\frac{1}{4\sqrt{3}} \sqrt{59 - \sqrt{2713}} \leq r \leq \frac{1}{\sqrt{3}}
$$ the following sharp inequality holds
$$
(1-|b_1|^2)\sum_{k=1}^\infty k|b_k|^2 r^{2k} \leq \frac{27}{8} r^4.
$$
\end{theorem}

Proof. At first let us remark it is enough to prove the theorem for the cases $r=R$ and $r=\sqrt{1/3}$. Indeed, for the analytic function
$$
\Psi(z)=(1-|b_1|^2)\sum_{k=1}^\infty k|b_k|^2 z^{2k}
$$ we have
$$
|\Psi(z)/z^4| \leq \frac{27}{8}
$$ for $r=R$ and $r=\sqrt{1/3}$ and therefore the inequality holds inside the ring $R \leq r \leq \sqrt{1/3}$.

It remains to prove the Theorem 5 in the case $r=R$.

We set $r=R$ and  consider three cases (we use these there cases due to technical reasons only, probably there exists a shorter proof).

{\bf Case 1:} $a \leq 3/5$. Remember that $a=\frac{3\sqrt{3}}{2}x(1-x^2)$. From here we see that $x \leq 0.287$ so that $(1/\sqrt{3}-x)/(1-x/\sqrt{3}) \ge 0.38 \ge R$  and we can apply Theorem 1.
In view of Theorem 1 we shall show that
$$(1 - a^2) \frac{27 R^2 (1 - x^2)^2 (2 x^2 +
   R^2 (1 + 2 (-3 + R^2) x^2 + x^4))}{8 (1- R^2 x^2)^4} - (27/8)R^4 \leq 0
$$ for $x \leq 1/4$.
We have
$$(1 - a^2) \frac{27 R^2 (1 - x^2)^2 (2 x^2 +
   R^2 (1 + 2 (-3 + R^2) x^2 + x^4))}{8 (1- R^2 x^2)^4} - (27/8)R^4 =
$$
$$
- \frac{27 R^2 x^4 }{32 (1 - R^2 x^2)^4}( 24.5695 \ldots - 103.49 \ldots x^2 + 159.036 \ldots x^4 -
$$
$$
- 99.9262 \ldots x^6 + 16.2316 \ldots x^{8} +
 3.88886 \ldots x^{10}).
$$ From here we see that the positive coefficients are dominating for $x \leq 1/4$ and consequently this expression is negative.

{\bf Case 2:}  $3/5  <a \leq 3/4$. In this case from Corollary 1 it follows that
$$
\sum_{k=2}^{\infty}k|b_k|^2r^{2k}\leq \frac{3(9-4a^2)^2 }{64a^4}\frac{16a^4}{9}\left(\log\frac{1}{1-r^2}-r^2\right)=
$$
$$
 \frac{(9-4a^2)^2 }{12}\left(\log\frac{1}{1-r^2}-r^2\right).
$$ Consequently
$$
(1-a^2)\sum_{k=1}^{\infty}k|b_k|^2r^{2k} \leq (1-a^2) \left(a^2r^2+\frac{(9-4a^2)^2 }{12}\left(\log\frac{1}{1-r^2}-r^2\right)\right).
$$
A little calculus shows that  maximum of the last expression is attained at the point $a=3/4$ and it is less than $27r^4/8$.

{\bf Case 3:}  $ a \ge 3/4$.

According to Corollary 2 the following inequality holds
\begin{equation}\label{B3}
\sum_{k=2}^\infty k|b_k|^2 r^{2k} \leq \frac{27}{8} r^4.
\end{equation}

From (\ref{B3}) it follows that
\begin{equation}\label{B5}
(1-|b_1|^2)\sum_{k=1}^\infty k|b_k|^2 r^{2k} \leq \max_{3/4 \leq a \leq 1}(1-a^2)\left(a^2 r^2+\frac{27}{8} r^4\right) \leq \frac{27}{8} r^4.
\end{equation}

{\bf Remark 5.} Routine and straightforward calculations show that the number $R$ in Theorem \ref{Th5} cannot be improved.\\

{\bf Problem 4.} How can one get inequalities analogous to the above ones in the intervals $[r_n,r_{n+1}], n\geq 1$ ?

{\bf Remark 6.} Computer experiments suggest that the inequality (\ref{Additional1}) can be replaced by $r \leq \sqrt{1/3}$.
If it is so then simple but routine calculations show that the lower bound $r=\sqrt{4/15}=0.5163...$ from Theorem \ref{Th2} can be replaced by the sharp
number $r=\sqrt{\rho}=0.39466...$
where $\rho$ is the positive root of the equation
$$
-4 - y + 81 y^2 + 642 y^3 - 564 y^4 + 1188 y^5 - 82 y^6 -
   5809 y^7 + 4581 y^8=0.
$$

  \subsection*{Acknowledgements}

  We thank the referee for his careful reading of our paper and his proposals that helped to ameliorate it.
The research of I. Kayumov was funded by the subsidy allocated to Kazan Federal University for the state assignment in the sphere of scientific activities, project No. 1.12878.2018/12.1.

I. R. Kayumov \\
Kazan Federal University \\
 Kremlevskaya 18 \\
  420 008 Kazan, Russia\\
e-mail: ikayumov@kpfu.ru

Karl-Joachim Wirths\\
Institut f\"ur Analysis\\
TU Braunschweig\\
38106 Braunschweig\\
Germany\\
e-mail: kjwirths@tu-bs.de
\end{document}